\documentclass[12pt]{article}

\usepackage{amsmath}
\usepackage{amssymb}

\newtheorem{theorem}{Theorem}
\newtheorem{lemma}{Lemma}

\newtheorem{definition}{Definition}

\def\Fp{\mathbb{F}_p}
\def\Fq{\mathbb{F}_q}

\begin{document}
\title{Additive properties of product sets in an arbitrary finite field\thanks{This paper was
supported by National Science Foundation grant under agreement No.
DMS-0635607.}.}
\author{Alexey Glibichuk}
\date{}
\maketitle
\begin{abstract} It is proved that for any two subsets $A$ and $B$
of an arbitrary finite field $\Fq$ such that $|A||B|>q$ the identity
$16AB=\Fq$ holds. Moreover, it is established that for every subsets
$X, Y\subset \Fq$ with the property $|X||Y|\geqslant 2q$ the
equality $8XY=\Fq$ holds.
\end{abstract}

\section{Introduction.}

Let $p$ be a prime, $m$ be a natural number, $\Fq$ be the finite
field of order $q=p^m$, and $\Fq^{*}$ be the multiplicative group of
$\Fq$, so that $\Fq^*=\Fq\setminus\{0\}$. For sets $X\subset\Fq$,
$Y\subset\Fq$, and for a (possibly, partial) binary operation
$*:\Fq\times\Fq\to\Fq$ we let
$$X*Y=\{x*y:\,x\in X,y\in Y\}.$$
We will write $XY$ instead of $X*Y$ if $*$ is multiplication in the
field; and, for an element $\lambda\in\Fq$, we write
$$\lambda\ast A=\{\lambda\}A$$
$$-A=(-1)\ast A=\{-a:a\in A\}.$$
For a set $X\subset\Fq$ and $k\in\mathbb{N}$ let
$$kX=\{x_1+\dots+x_k:\,x_1,\dots,x_k\in X\},$$
$$X^k=\{x_1\dots x_k:\,x_1,\dots,x_k\in X\}.$$
Let also denote the cardinality of the given set $X$ as $|X|$. For
given natural numbers $N$ the notation $NXY$ should be understood as
$N$-fold sum of the product set $XY$. Let us consider the following
definitions.
\begin{definition} The set $X$ is said to be \textbf{symmetric} if
$X=-X$.
\end{definition}
\begin{definition} The set $X$ is said to be \textbf{antisymmetric} if
$X\cap(-X)=\emptyset.$
\end{definition}

A set $A$ is called an (additive) basis  of order $k$ (for $\Fq$) if
$kA=\Fq$. Observe that any basis of order $k$ is also a basis of any
order $k'>k$. The general problem that will be discussed in this
paper is whether, for given integers $t<q, N$ and two sets $A$ and
$B$, the set $AB$ is a basis of order $N$ if $|A||B|\geq t$?

The first machinery, allowing one to prove sum-product results on
finite fields was developed in the paper of J.~Bourgain, N.~Katz and
T.~Tao(\cite{BKT}).

The author of this paper proved the following two
statements(\cite{Gl}, Theorems 1 and 2).

\begin{theorem}\label{first8AB1}Let $A$ and $B$ be subsets of the
field $\Fp$ for some prime $p$. If the set $B$ is antisymmetric and
$|A||B|>p$ then $8AB=\Fp.$
\end{theorem}
\begin{theorem}\label{first8AB2}Let $A$ and $B$ be subsets of the
field $\Fp$ for some prime $p$. If the set $B$ is symmetric and
$|A||B|>p$ then $8AB=\Fp.$
\end{theorem}

In the joint paper with S.V. Konyagin(\cite{GK}, Lemmas 2.1 and 2.2
)we established the following two results.

\begin{theorem}\label{first8ABarb} If $A\subset\Fp$, $B\subset\Fp$
for some prime $p$, and
$|A|\cdot\lceil|B|/2\rceil>p$ then $8AB=\Fp$.
\end{theorem}
\begin{theorem}\label{first16AB} If $A\subset\Fp$, $B\subset\Fp$
for some prime $p$, and $|A||B|>p$ then $16AB=\Fp$.
\end{theorem}

In this paper extensions of Theorems \ref{first8AB1}-\ref{first16AB}
will be obtained. We shall establish the following four theorems.

\textbf{Theorem \ref{thm:Bisantysymm}} \emph{If $A\subset\Fq$ and
$B\subset\Fq$ are such that $B$ is antisymmetric and $|A||B|>q$ then
$8AB=\Fq.$}

\textbf{Theorem \ref{thm:Bissumm}} \emph{Assume that $A\subset\Fq$
and $B\subset\Fq$ are such that $B$ is symmetric. If also $|A||B|>q$
then $8AB=\Fq.$}

\textbf{Theorem \ref{main}} \emph{Let $A, B\subset\Fq$ be arbitrary
subsets with $|A||B|>q.$ Then we have $16AB=\Fq.$}

\textbf{Theorem \ref{main1}} \emph{Let $A,B\subset\Fq$ be arbitrary
subsets with $|A||B|\geqslant 2q.$ Then we have $8AB=\Fq.$}

Constant $16$ in the Theorem \ref{main} is most likely not best
possible, it is demonstrated by Theorem \ref{main1} and recent
result of D. Hart and A. Iosevich(\cite{IH}). They established that
\begin{theorem}\label{HI}For every subset $A \subset\Fq$
such that $|A| \ge Cq^{\frac{1}{2}+\frac{1}{2d}}$ for $C$
sufficiently large the identity $dA^2=\Fq^{*}$ holds.
\end{theorem}

Applying Theorem \ref{HI} with $d=1$ we see that the constant in the
Theorem \ref{main} can be significantly improved when $A=B$ and
$|A|>Cq^{\frac{3}{4}}.$ D. Hart and A. Iosevich in the same paper
have conjectured that if
$|A|>C_{\varepsilon}q^{\frac{1}{2}+\varepsilon}$ for some constant
$C_{\varepsilon}$ and $\varepsilon>0$ then $2A^2=\Fq.$ However,
condition $|A||B|>q$ in the Theorem \ref{main} is sharp. Indeed, if
$|A||B|=q$ then result similar to the Theorem \ref{main} cannot
hold. It is sufficient to consider sets $A=\Fq, B=\{0\}$ or make
$A=B$ to be a subfield of order $\sqrt{q}$ when $q=p^m$ and $m$ is
even, to verify this statement. To construct a less trivial
counterexample let us consider two natural numbers $k$ and $l$ such
that $k+l=m.$ Let us take a primitive element $\xi\in\Fq^{*}$ and
consider sets
$$A=\{x_0+x_1\xi+\ldots+x_{k-1}\xi^{k-1}:(x_0,\ldots,x_{k-1})
\in\underbrace{\Fp\times\ldots\times\Fp}_{k}:=\Fp^{k}\},$$
$$B=\{x_0+x_1\xi+\ldots+x_{l-1}\xi^{l-1}:(x_0,x_1,\ldots,x_{l-1})\in\Fp^{l}\}$$
and
$$C=\{x_0+x_1\xi+\ldots+x_{m-2}\xi^{m-2}:(x_0,x_1,\ldots,x_{m-2})\in\Fp^{m-1}\}$$
where $\Fp\subset\Fq$ is a subfield of $\Fq$ of cardinality $p$.
Then one can obviously observe that $|A||B|=q,$ $AB\subset C\neq\Fq$
and $C$ is closed under addition.

\section{Preliminary results.}

Lemmas \ref{technical}, \ref{required}, \ref{IABlemma} are
extensions of Lemmas 1, 2, 3 from \cite{Gl}. Their proofs are due to
arguments used in corresponding lemmas.

\begin{lemma}\label{technical}
Let $A\subset\Fq$, $B\subset\Fq$ be arbitrary nonempty subsets. Then
there is an element $\xi\in\Fq^{*}$ such that
\begin{equation}\label{eqn:AplusxiBislarge}
|A+\xi B|\geqslant\frac{|A||B|(q-1)}{|A||B|-(|A|+|B|)+q}
\end{equation}
and
\begin{equation}\label{eqn:AminusxiBislarge}
|A-\xi B|\geqslant\frac{|A||B|(q-1)}{|A||B|-(|A|+|B|)+q}.
\end{equation}
\end{lemma}
\textbf{Proof.} Let us consider an arbitrary elements
$\xi\in\Fq^{*}$ and $s\in\Fq.$ Denote
$$f_{\xi}^{+}(s)=|\{(a,b)\in A\times B:a+b\xi=s\}|,$$
$$f_{\xi}^{-}(s)=|\{(a,b)\in A\times B:a-b\xi=s\}|.$$
It obviously follows that
$$\sum_{s\in
\Fq}(f_{\xi}^{+}(s))^{2}=|\{(a_1,b_1,a_2,b_2)\in A\times B\times
A\times B:a_1+b_1\xi = a_2+b_2\xi \}|$$
$$=|A||B|+|\{(a_1,b_1,a_2,b_2)\in A\times B\times A\times B:a_1
\neq a_2,a_1+b_1\xi=a_2+b_2\xi \}|,
$$
$$\sum_{s\in\Fq}(f_{\xi}^{-
}(s))^{2}=|\{(a_1,b_1,a_2,b_2)\in A\times B\times A\times B:
a_1-b_1\xi = a_2-b_2\xi \}|$$
$$=|A||B|+|\{(a_1,b_1,a_2,b_2)\in A\times B\times A\times B: a_1\neq
a_2,a_1-b_1\xi=a_2-b_2\xi \}|.$$ Therefore, $\sum_{s\in
\Fq}(f_{\xi}^{+}(s))^{2}=\sum_{s\in\Fq}(f_{\xi}^{- }(s))^{2},$ and
it is enough to consider only sum with values of $f_{\xi}^{+}(s).$
It is easy to see that for every $a_1,a_2\in A$ and $b_1,b_2\in B$
with $a_1\neq a_2$ there is only one element $\eta\neq 0$ such that
$a_1+b_1\eta=a_2+b_2\eta.$ Thus,
$$\sum_{\xi\in\Fq^{*}}\sum_{s\in\Fq}(f_{\xi}^{+}(s))^{2}=
|A||B|(q-1)+|A||B|(|A|-1)(|B|-1).$$ Therefore, there is
$\xi\in\Fq^{*}$ such that
\begin{equation}\label{eqn:sum_sfplus}
\sum_{s\in\Fq}(f_{\xi}^{+}(s))^{2}\le
|A||B|+\frac{|A||B|(|A|-1)(|B|-1)}{q-1}.
\end{equation}
By Cauchy-Schwartz
\begin{equation}\label{eqn:causchwcoroll1}
\left(\sum_{s\in\Fq}f_{\xi}^{+}(s)\right)^{2}\le |A+\xi
B|\sum_{s\in\Fq}(f_{\xi}^{+}(s))^{2},
\end{equation}
\begin{equation}\label{eqn:causchwcoroll2}
\left(\sum_{s\in\Fq}f_{\xi}^{-}(s)\right)^{2}\le |A-\xi
B|\sum_{s\in\Fq}(f_{\xi}^{-}(s))^{2}.
\end{equation}
Moreover, it obviously follows that
$$\sum_{s\in\Fq}f_{\xi}^{+}(s)=|A||B|,$$
$$\sum_{s\in\Fq}f_{\xi}^{-}(s)=|A||B|.$$
Now from (\ref{eqn:sum_sfplus}), (\ref{eqn:causchwcoroll1}) and
(\ref{eqn:causchwcoroll2}) one can deduce a desired inequality:
$$|A+\xi B|\ge
\frac{|A|^2|B|^2}{|A||B|+\frac{|A||B|(|A|-1)(|B|-1)}{q-1}}=\frac{|A||B|(q-1)}{
|A||B|-(|A|+|B|)+q}
$$
and
$$|A-\xi B|\ge\frac{|A||B|(q-1)}{|A||B|-(|A|+|B|)+q}.$$
Lemma 1 is proved. $\blacksquare$

\begin{lemma}\label{required} Let $A$ and $B$ be subsets
of field $\Fq$ with $|A||B|>q$. Then there is $\xi\in\Fq^{*}$ such
that
\begin{equation}\label{eqn:AplusxiBqhalf}
|A+\xi B|>\frac{q}{2}
\end{equation}
and
\begin{equation}\label{eqn:AminusxiBqhalf}
|A-\xi B|>\frac{q}{2}.
\end{equation}
\end{lemma}
\textbf{Proof.} Let us apply Lemma \ref{technical}. It states that
there is $\xi\in\Fq^{*}$ such that (\ref{eqn:AplusxiBislarge}) and
(\ref{eqn:AminusxiBislarge}) hold. Clearly, we have
$$
\frac{|A||B|(q-1)}{|A||B|-(|A|+|B|)+q}\geqslant\frac{|A||B|(q-1)}{|A||B|+(q-2)}.
$$
Let us consider the difference
$$s=\frac{|A||B|(q-1)}{|A||B|+(q-2)}-\frac{q}{2}=\frac{(q-2)(|A||B|-q)}{2(|A||B|+(q-2))}.$$
It is clear that $s>0$ when $|A||B|>q$ and $q\neq 2.$ If $q=2$ then
the condition $|A||B|>q$ implies that at least one of the subsets
$A$ or $B $ is equal to $\Fq.$ Lemma \ref{required} is proved.
$\blacksquare$

\begin{definition} For two subsets $A\subset \Fq, B\subset\Fq$ denote
$$I(A,B)=\{(b_1-b_2)\cdot a_1+(a_2-a_3)\cdot b_3:a_1,a_2,a_3\in A, b_1,b_2,b_3\in B\}.$$
\end{definition}

\begin{lemma}\label{IABlemma} Consider two subsets $A\subset\Fq$ and
$B\subset\Fq$. If for some $\xi\in\Fq^{*}$
$$|A+\xi B|<|A||B|$$
then $$|I(A,B)|\geqslant |A+\xi B|.$$
\end{lemma}
\textbf{Proof.} If $|A+\xi B|<|A||B|$ then there are elements $a_1,
a_2\in A$ and $b_1, b_2\in B$ such that $(a_1,b_1)\neq (a_2,b_2)$
and
\begin{equation}\label{eqn:nontriverelation}
(a_1-a_2)+(b_1-b_2)\cdot\xi=0.
\end{equation}
It is clear that $b_1\neq b_2$. Let us consider the set
$$S=(b_1-b_2)\cdot (A+\xi B)=\{(b_1-b_2)\cdot b:b\in A+\xi B\}.$$
It is obviously follows that $|S|=|A+\xi B|$ and every element of
$S$ can be rewritten in the form $$s=(b_1-b_2)\cdot a+(b_1-b_2)\cdot
b\xi$$ with $a\in A$ and $b\in B$. From (\ref{eqn:nontriverelation})
one can easily deduce that
$$s=(b_1-b_2)\cdot a+(a_2-a_1)\cdot b.$$
Therefore, $S\subset I(A,B)$ and lemma follows. $\blacksquare$

\begin{lemma}\label{Abigsumbig} Assume that $X\subset\Fq$ with
$|X|>\frac{q}{2},$ then $X+X=\Fq.$
\end{lemma}
\textbf{Proof.} Let us take an arbitrary element $x\in\Fq$ and
consider a set $x-X.$ From $|x-X|=|X|>\frac{q}{2}$ one can obviously
prove that sets $x-X$ and $X$ have nonempty intersection, so there
are elements $x_1, x_2\in X$ such that $x-x_1=x_2\Leftrightarrow
x=x_1+x_2.$ Lemma now follows. $\blacksquare$

\begin{lemma}\label{AplusAhaslargeregsub} Let $A$ be any subset of
$\Fq.$ If $|A|\not\equiv 2\pmod 3$ then there is a symmetric or
antisymmetric subset $S\subset A$ with $|S|\geqslant
\frac{2}{3}|A|$. If $|A|\equiv 2\pmod 3$ then one can find either
symmetric or antisymmetric subset $S\subset A$ with $|S|\geqslant
\frac{2}{3}|A|-\frac{1}{3}.$
\end{lemma}
\textbf{Proof.} Let us define a set $A_1=\{x\in A:-x\notin A\}$. It
is an antisymmetric subset of $A$. Consider a set of subsets
$\mathcal{S}=\{\{a_1,a_2\}:a_1\in A, a_2\in A, a_1=-a_2\}.$ It is
clear that one can choose one element from each of the sets from
$\mathcal{S}$ and form a new set $A_2$ from those elements. It is
easy to observe that $A_2\cap A_1=\emptyset,$ $0\in A_2$ if $0\in A$
and $A_2\setminus\{0\}$ is antisymmetric. Let us define a subset
$A_3=A\setminus(A_1\sqcup A_2).$ It is an antisymmetric subset of
$A$ with cardinality $|A_3|=|A_2|-1$ if $0\in A$ and $|A_3|=|A_2|$
otherwise, such that $0\notin A_3$ and $A_2\sqcup A_3$ is the
maximal symmetric subset of $A$. We have split the set $A$ into
three nonintersecting parts: $A=A_1\sqcup A_2\sqcup A_3$.

If $|A_1|<\frac{1}{3}|A|$ then $|A_2\sqcup
A_3|\geqslant\frac{2}{3}|A|$ and Lemma \ref{AplusAhaslargeregsub}
follows with symmetric $S=A_2\sqcup A_3$.

If $0\notin A$ and $|A_1|\geqslant\frac{1}{3}|A|$ then $|A_2|=|A_3|$
and $|A_3|<\frac{1}{3}|A|$. Assuming $S$ to be an antisymmetric
subset $A_1\sqcup A_2$ we complete the proof of Lemma
\ref{AplusAhaslargeregsub}.

Assume that $|A_1|\geqslant\frac{1}{3}|A|$ and $0\in A$. If
$|A_3|\geqslant\frac{1}{3}|A|$ then $|A_1\sqcup
A_3|\geqslant\frac{2}{3}|A|$ and Lemma \ref{AplusAhaslargeregsub} is
proved by letting $S$ to be antisymmetric subset $A_1\sqcup A_3$.

It is left to prove Lemma \ref{AplusAhaslargeregsub} when
\begin{equation}\label{eqn:A1isbig}
|A_1|\geqslant\frac{1}{3}|A|,
\end{equation}
\begin{equation}\label{eqn:A3issmall}
|A_3|<\frac{1}{3}|A|
\end{equation}
and $0\in A.$ Let us consider three cases.

Case 1. $|A|=3k$ for some natural $k$. Taking into account
(\ref{eqn:A1isbig}) and (\ref{eqn:A3issmall}) one can see that
$|A_3|\leqslant k-1$ and therefore $|A_1\sqcup A_2|\geqslant 2k+1.$
By defining $S=(A_1\sqcup A_2)\setminus\{0\}$ ($S$ is antisymmetric)
we complete the proof of Lemma \ref{AplusAhaslargeregsub}.

Case 2. $|A|=3k+1$ for some natural $k$. Again, using
(\ref{eqn:A1isbig}) and (\ref{eqn:A3issmall}) one can deduce that
$|A_3|\leqslant k$. If $|A_3|\leqslant k-1$ then assuming
$S=(A_1\sqcup A_2)\setminus\{0\}$ we get a required antisymmetric
subset. If $|A_3|=k$ then $|A_2|=k+1$ and $|A_1|=k.$ Note that the
identity $|A_1|=k$ contradicts inequality (\ref{eqn:A1isbig}). We
are done.

Case 3. $|A|=3k+2$ for some natural $k$. Using (\ref{eqn:A1isbig})
and (\ref{eqn:A3issmall}) one can easily deduce that $|A_3|\leqslant
k$ and $|A_1|\geqslant k+1$. If $|A_3|\leqslant k-1$ then
$|A_1\sqcup A_2|=|A\setminus A_3|\geqslant 2k+3$. Letting $S$ to be
an antisymmetric subset $(A_1\sqcup A_2)\setminus\{0\}$ we observe
that $|S|\geqslant 2k+2>\frac{1}{3}|A|$ and we are done with better
bound on $|S|$. In case when $|A_3|=k$ it is easy to see that
$|A_2|=k+1$ and $|A_1|=k+1$. Assuming $S$ to be a symmetric subset
$A_2\sqcup A_3$ we complete the proof of Lemma
\ref{AplusAhaslargeregsub}. $\blacksquare$

\begin{definition} For every subset $X\subset\Fq$ its
\textbf{symmetry group} (it is denoted as $Sym_1(X)$) is defined by
the identity
$$Sym_1(X)=\{h:\{h\}+X=X\}.$$
\end{definition}

We shall use the following theorem (see \cite{TV}, theorem 5.5 or
\cite{Kn}).
\begin{theorem}\label{Kneser}(Kneser) For every subsets $X,
Y\subset\Fq$ we have
$$|X+Y|\geqslant |X+Sym_1(X+Y)|+|Y+Sym_1(X+Y)|-|Sym_1(X+Y)|\geqslant$$
$$\geqslant |X|+|Y|-|Sym_1(X+Y)|.$$
\end{theorem}
\begin{lemma}\label{symetrystruct} Given a subset $X\subset\Fq.$ Let us
take any subgroup $G$ of the group $Sym_1(X).$ Then $X$ is a union
of additive cosets of $G.$
\end{lemma}
\textbf{Proof.} One can easily observe that $Sym_1(X)$ is an
additive subgroup. It is sufficient to prove that every coset of the
subgroup $G$ either is a subset of $X$ or has an empty intersection
with $X$. Suppose that some coset $x+G$ has nonempty intersection
with $X$. Let us take an arbitrary element $y\in X\cap(x+G)$. By
definition of $y$ a coset $y+G=x+G,$ but from symmetry of $X$ it
follows that $y+G\subset X.$ Lemma \ref{symetrystruct} is proved.
$\blacksquare$

\begin{lemma}\label{AplusAislarge} Let $B$ be an arbitrary subset of
$\Fq$ such that $|B|\geqslant 2.$ Then one of the following two
alternatives holds
\begin{itemize}
\item[(i)] $|B+B|\geqslant\frac{3}{2}|B|,$
\item[(ii)] there is an additive subgroup $G\subset\Fq$ such that $B\subset
b+G$ for some $b\in B$ and $|B|>\frac{2}{3}|G|$. Moreover, in this
case $B+B=2b+G.$
\end{itemize}
\end{lemma}
\textbf{Proof.} Application of Theorem \ref{Kneser} for sets $X=Y=B$
implies
\begin{equation}\label{eqn:kneserapp}
\begin{split}
|B+B|\geqslant 2&|B+Sym_1(B+B)|-|Sym_1(B+B)|\geqslant\\
&2|B|-|Sym_1(B+B)|.
\end{split}
\end{equation}

Since $Sym_1(B+B)$ is an additive subgroup of $\Fq$ then there is an
integer $0\leqslant l\leqslant n$ such that $|Sym_1(B+B)|=p^l.$
Observe that $Sym_1(B+B)\subset Sym_1(B+Sym_1(B+B)).$ Now from Lemma
\ref{symetrystruct} clearly follows that $|B+Sym_1(B+B)|=mp^l$ for
some natural $m$. Again, using (\ref{eqn:kneserapp}) we can see that
\begin{equation}\label{eqn:AplusAwithm}
|B+B|\geqslant(2m-1)p^l.
\end{equation}

Assume that the inequality $|B+B|<\frac{3}{2}|B|$ holds. Then we
deduce from (\ref{eqn:kneserapp}) that
$$\frac{3}{2}|B|>2|B|-|Sym_1(B+B)|\Leftrightarrow |B|<2p^l$$
and therefore $|B+B|<\frac{3}{2}\cdot 2p^l=3p^l.$ Combining the last
inequality with (\ref{eqn:AplusAwithm}) we obtain the condition
$2m-1<3$ and therefore $m$ can take on one value: $m=1.$ When $m=1$
one can observe that $|B+Sym_1(B+B)|=|Sym_1(B+B)|=p^l$. Take an
arbitrary element $b\in B$ and consider the set $B^{'}=B-b$. It is
clear, that $B^{'}+Sym_1(B+B)=Sym_1(B+B)$ and therefore
$B^{'}\subset Sym_1(B+B).$ Recalling definition of $B^{'}$ we obtain
a relation $B\subset b+Sym_1(B+B).$ By (\ref{eqn:AplusAwithm}) one
can deduce the inequality $|B+B|\geqslant p^l.$ Observing that
$B+B\subset 2b+Sym_1(B+B)$ we can obtain the relation
$|B+B|=|Sym_1(B+B)|=p^l.$ Now it is clear that if
$|B|\leqslant\frac{2}{3}p^l$ then the inequality
$|B+B|\geqslant\frac{3}{2}|B|$ holds, otherwise we get the
alternative $(ii).$ To finish the proof of the Lemma
\ref{AplusAislarge} we need to observe that according to Lemma
\ref{Abigsumbig} $B+B=2b+Sym_1(B+B)$ when $|B|>\frac{2}{3}p^l.$
Lemma \ref{AplusAislarge} now follows. $\blacksquare$

\section{Proofs of theorems \ref{thm:Bisantysymm}-\ref{main1}.}

\begin{theorem}\label{thm:Bisantysymm} If $A\subset\Fq$ and
$B\subset\Fq$ are such that $B$ is antisymmetric and $|A||B|>q$ then
$8AB=\Fq.$
\end{theorem}
\textbf{Proof.} Let us apply Lemma \ref{required}. It states that
there is an element $\xi\in\Fq^{*}$ such that
(\ref{eqn:AplusxiBqhalf}) and (\ref{eqn:AminusxiBqhalf}) hold. From
(\ref{eqn:AplusxiBqhalf}) one can easily derive that $(A+\xi b)\cap
(-A-\xi B)\neq\emptyset$ and, therefore, there are elements $a_1,
a_2\in A, b_1, b_2\in B$ with $a_1+b_1\xi=-(a_2+b_2\xi)$. Thus,
\begin{equation}\label{eqn:localexprforxi}
\xi=-\frac{a_1+a_2}{b_1+b_2}.
\end{equation}
The expression (\ref{eqn:localexprforxi}) is correct because $B\cap
(-B)=\emptyset$ and denominator of the fraction in this formula is
not equal to zero. From (\ref{eqn:AminusxiBqhalf}) it follows that
$$\left|\left\{a_3+\frac{a_1+a_2}{b_1+b_2}b_3:a_3\in A,b_3\in B\right\}\right|>\frac{q}{2}\Leftrightarrow$$
$$|\{a_3(b_1+b_2)+b_3(a_1+a_2):a_3\in A,b_3\in B\}|>\frac{q}{2}.$$
Therefore, $|4AB|>\frac{q}{2}$ and Lemma \ref{Abigsumbig} gives us
the desired statement. $\blacksquare$

\begin{theorem}\label{thm:Bissumm} Assume that $A\subset\Fq$ and
$B\subset\Fq$ are such that $B$ is symmetric. If also $|A||B|>q$
then $8AB=\Fq.$
\end{theorem}
\textbf{Proof.} Applying Lemma \ref{required} one can find an
element $\xi\in\Fq^{*}$ such that $|A+\xi B|>\frac{q}{2}$. Moreover,
from restrictions on sets $A$ and $B$ one can see that $|A+\xi
B|\leqslant q<|A||B|$ and we can apply Lemma \ref{IABlemma} that
gives us the following:
$$|I(A,B)|\geqslant |A+\xi B|>\frac{q}{2}.$$
Taking into account that $B=-B$ one can derive that $I(A,B)\subset
4AB$ and $|4AB|>\frac{q}{2}.$ Now Theorem \ref{thm:Bissumm} follows
from Lemma \ref{Abigsumbig}. $\blacksquare$

\begin{theorem}\label{main} Let $A, B\subset\Fq$ be arbitrary
subsets with $|A||B|>q.$ Then we have $16AB=\Fq.$
\end{theorem}
\textbf{Proof.} Let us apply Lemma \ref{AplusAislarge} for the set
$B$. If $(ii)$ holds then $B+B=2b+G$ for some $b\in B$ and an
additive subgroup $G\subset\Fq$. It is easy to see that every coset
of an additive subgroup is an antisymmetric or a symmetric subset.
Then application of Theorem \ref{thm:Bisantysymm} or Theorem
\ref{thm:Bissumm} for sets $A$ and $B+B$ gives us Theorem
\ref{main}.

Assume now that
\begin{equation}\label{eqn:alerntwo}
|B+B|\geqslant\frac{3}{2}|B|
\end{equation}
i. e. alternative $(i)$ holds. If $|B+B|\not\equiv 2\pmod 3$ then
application of Lemma \ref{AplusAhaslargeregsub} gives us a subset
$S\subset B+B$ such that $|S|\geqslant\frac{2}{3}|B+B|$ and $S$ is
either symmetric or antisymmetric. By (\ref{eqn:alerntwo}) we
observe that $|S|\geqslant |B|$. Application of Theorem
\ref{thm:Bisantysymm} or Theorem \ref{thm:Bissumm} for sets $A$ and
$S$ allows one deduce Theorem \ref{main}.

It is left to consider the case when $|B+B|\equiv 2\pmod 3$ and the
inequality (\ref{eqn:alerntwo}) holds. Assume that $|B+B|=3k+2$ for
some natural $k$. Lemma \ref{AplusAhaslargeregsub} states that there
is either symmetric or antisymmetric subset $S\subset B+B$ with
$|S|\geqslant\frac{2}{3}|B+B|-\frac{1}{3}=2k+1.$ Moreover, by
(\ref{eqn:alerntwo}) we can deduce that $2k+1\geqslant
|B|-\frac{1}{3}$ and, therefore $|B|\leqslant 2k+1$. Now it is easy
to see that $|S|\geqslant 2k+1\geqslant |B|.$ Using Theorem
\ref{thm:Bisantysymm} or Theorem \ref{thm:Bissumm} for sets $A$ and
$S$ we complete the proof of Theorem \ref{main}. $\blacksquare$

\begin{theorem}\label{main1} Let $A,B\subset\Fq$ be arbitrary
subsets with $|A||B|\geqslant 2q.$ Then we have $8AB=\Fq.$
\end{theorem}
\textbf{Proof.} Our aim is to extract from the set $B$ a
sufficiently large symmetric or antisymmetric subset $S\subset B$.
Lemma \ref{AplusAhaslargeregsub} states that there is a symmetric or
antisymmetric subset $S\subset B$ with $|S|\geqslant
\frac{2}{3}|B|-\frac{1}{3}$. Let us notice that the equality $|B|=2$
holds when $|A||B|=2q$ and, therefore, $A=\Fq$, so we can assume
that $|B|>2.$ Observe that in this case
$\frac{2}{3}|B|-\frac{1}{3}>\frac{1}{2}|B|$ and we have $|A||S|>
\frac{1}{2}|A||B|\geqslant q$ and Theorem \ref{main1} now follows
from Theorems \ref{thm:Bisantysymm} and \ref{thm:Bissumm}.
$\blacksquare$


\begin{thebibliography}{}
\bibitem{BKT} J.~Bourgain, N.~Katz, T.~Tao, \emph{A sum-product estimate in
finite fields and their applications}, Geom and Funct. Anal., {\bf
14} (2004), 27--57.

\bibitem{Gl} A. A. Glibichuk, \emph{Combinational properties of sets of residues modulo a prime
and the Erd\H{o}s-Graham problem,} Mat. Zametki, {\bf 79} (2006),
384--395; translation in: Math. Notes, {\bf 79} (2006), 356--365.

\bibitem{GK} A. A. Glibichuk, S.V. Konyagin, \emph{Additive properties of product sets
in fields of prime order}, Centre de Recherches Math\'{e}matiques
Proceedings and Lecture Notes, vol. 43, pp. 279-286.

\bibitem{IH} D.Hart, A. Iosevich, \emph{Sums and products in finite fields:
an integral geometric viewpoint,}
preprint.

\bibitem{Kn} M. Kneser, \emph{Absch\"{a}tzungen der asymptotischen Dichte
von Summenmengen}, Math. Z, vol. 58, 1953, pp. 459--484.

\bibitem{TV} T.~Tao, V.~Vu, \emph{Additive combinatorics},
Cambridge Univ.~Press, Cambridge, 2006.

\end{thebibliography}
\end{document}